\documentclass[11pt,twoside]{article}
\title{} \author{} \date{}
\usepackage{latexsym,amssymb,times}
\input amssym.def
\newtheorem{te}{Theorem}[section]
\newtheorem{prop}[te]{Proposition}

\newtheorem{fac}[te]{Fact}

\newtheorem{ex}[te]{Example}

\def\dok{\noindent{\bf Proof. }}
\def\kdok{\hfill $\Box$ \par \vspace*{2mm} }
\def\a{\alpha}
\def\b{\beta}
\def\k{\kappa}
\def\l{\lambda}

\def\o{\omega}

\def\c{\mathfrak{c}}
\def\N{{\mathbb N}}
\def\Q{{\mathbb Q}}
\def\R{{\mathbb R}}
\def\Z{{\mathbb Z}}

\def\BL{{\mathbb L}}
\def\BC{{\mathbb C}}

\def\C{{\mathcal C}}

\def\O{{\mathcal O}}
\def\F{{\mathcal F}}

\def\L{{\mathcal L}}
\def\U{{\mathcal U}}

\def\CS{{\mathcal S}}
\def\CX{{\mathcal X}}
\def\CY{{\mathcal Y}}
\def\CR{{\mathcal R}}
\def\CW{{\mathcal W}}

\def\la{\langle}
\def\ra{\rangle}

\def\Top{\mathop{{\mathrm{Top}}}\nolimits}

\def\otp{\mathop{{\mathrm{otp}}}\nolimits}

\def\Aut{\mathop{{\mathrm{Aut}}}\nolimits}

\def\Sym{\mathop{{\mathrm{Sym}}}\nolimits}
\def\Conv{\mathop{{\mathrm{Conv}}}\nolimits}
\begin{document}
\thispagestyle{plain}
\begin{center}
           {\large \bf WEAKLY AND STRONGLY REVERSIBLE SPACES}
\end{center}
\begin{center}
{\bf Milo\v s S.\ Kurili\'c}\footnote{Department of Mathematics and Informatics, Faculty of Sciences, University of Novi Sad,
                                      Trg Dositeja Obradovi\'ca 4, 21000 Novi Sad, Serbia.
                                      e-mail: milos@dmi.uns.ac.rs}
\end{center}
\begin{abstract}
\noindent
A topological space ${\mathcal X}$ is 
{\it reversible} iff each continuous bijection (condensation) $f: {\mathcal X} \rightarrow {\mathcal X}$ is a homeomorphism;
{\it weakly reversible} iff whenever ${\mathcal Y}$ is a space and there are condensations $f:{\mathcal X} \rightarrow {\mathcal Y}$ and $g:{\mathcal Y} \rightarrow {\mathcal X}$,
there is a homeomorphism  $h:{\mathcal X} \rightarrow {\mathcal Y}$;
{\it strongly reversible} iff each bijection $f: {\mathcal X} \rightarrow {\mathcal X}$ is a homeomorphism.
We show that the class of weakly reversible non-reversible spaces
is disjoint from the class of sequential spaces in which each sequence has at most one limit (containing e.g.\ metrizable spaces).
On the other hand, the class of strongly reversible topologies contains only discrete topologies,  antidiscrete topologies and natural generalizations of the cofinite topology.


\end{abstract}

\section{Preliminaries}
Reversibility is a topological property detected and initially investigated by Rajagopalan and Wilansky in \cite{Raj}:
a space $\CX =(X,\O)$ (or a topology $\O \in \Top _X$) is called {\it reversible}
iff each continuous bijection (condensation) $f: \CX \rightarrow \CX$ is a homeomorphism.
In this note we consider two natural variations of that property.
A space $\CX$ will be called {\it strongly reversible}
iff each bijection $f: \CX \rightarrow \CX$ is a homeomorphism;
$\CX$ will be called {\it weakly reversible}
iff whenever $\CY$ is a space and there are condensations $f:\CX \rightarrow \CY$ and $g:\CY \rightarrow \CX$,
there is a homeomorphism  $h:\CX \rightarrow \CY$.
Denoting the corresponding classes of spaces by $\CR$, $\CS\CR$ and $\CW\CR$ it is easy to check that $\CS\CR \subset \CR \subset \CW\CR$
and that the class $\CW\CR$ satisfies the Cantor-Schr\"{o}eder-Bernstein property for condensations.

Since the aforementioned topological properties are related to the existence of bijections between spaces,
w.l.o.g.\ we can restrict our attention to the complete lattice $\la \Top _X,\subset \ra$ of all topologies on a set $X$ (or, moreover, on the cardinal $\k=|X|$).
For topologies $\O ,\O ' \in \Top _X$ we will write $\O  \cong \O '$ iff there is a homeomorphism $f:(X,\O )\rightarrow (X,\O ')$,
which means that $f[\O ]=\O'$, where slightly abusing notation we define $f[\O ]:=\{ f[O]:O\in \O \}$.
Clearly, $[\O ]_{\cong }:=\{ \O '\in \Top _X : \O ' \cong \O\}=\{ f[\O ]: f\in \Sym _X\}$ and by \cite{Raj} we have
\begin{fac}\label{T107}
For a topology $\O \in \Top _X$ the following conditions are equivalent

(i) $\O$ is reversible,

(ii) There is no topology $\O ' \in [\O ]_{\cong}$ such that $\O ' \subsetneq \O $,

(iii) There is no topology $\O ' \in [\O ]_{\cong}$ such that $\O \subsetneq \O '$,

(iv) The poset $\la [\O ]_{\cong},\subset \ra$ is an antichain (different elements are incomparable).
\end{fac}
Investigating the diversity of non-homeomorphic topologies on a set $X$
it is natural to compare the topologies from $\Top _X$
saying that $\O_1 \preccurlyeq \O _2$ iff there is $\O \in \Top _X$ such that $\O _1 \cong \O \subset \O _2$.
The following characterization has a direct proof.
\begin{fac}\label{T103}
For topologies $\O _1,\O _2 \in \Top _X$ the following conditions are equivalent

(i)  $\O_1 \preccurlyeq \O _2$,

(ii) There is a topology $\O \in \Top _X$ such that $\O _1 \subset \O \cong \O _2$,

(iii) There is a condensation $f: (X, \O _2) \rightarrow (X,\O _1)$.
\end{fac}
So we obtain the {\it condensational preorder} $\la \Top _X ,\preccurlyeq \ra $
and, defining the {\it condensational equivalence} relation by $\O _1 \sim \O _2$ iff $\O _1 \preccurlyeq \O _2$ and $\O _2 \preccurlyeq \O _1$,\footnote{
such topologies are called {\it bijectively related} in \cite{Doy}}
we obtain its antisymmetric quotient, the {\it condensational order} $\la \Top _X \!/\!\sim ,\leq \ra$.
Here, for topologies $\O,\O _1,\O _2\in \Top _X$ we have $[\O ]_{\cong }\subset[\O ]_{\sim }:=\{ \O '\in \Top _X : \O ' \sim \O\}$
and $[\O _1]_{\sim} \leq [\O _2]_{\sim}$ iff $\O _1 \preccurlyeq \O _2$.

In order to characterize weakly reversible topologies,
for a set $\Omega \subset \Top _X$ by $\Conv (\Omega )$ we denote the minimal convex subset of the lattice $\la \Top _X ,\subset\ra$ containing $\Omega$;
namely, $\Conv (\Omega ):=\{ \O \in \Top _X : \exists \O _1,\O _2\in \Omega \;\;\O _1 \subset \O \subset \O _2\}$.
\begin{prop}\label{T104}
For each topology $\O \in \Top _X$ we have\\[-7mm]
\begin{itemize}\itemsep=-1.5mm
\item[\rm (a)] $[\O ]_\sim$ is a convex set in the lattice $\la \Top _X ,\subset\ra$;
\item[\rm (b)] $[\O ]_\sim =\Conv ([\O ]_{\cong})$;
\item[\rm (c)] The following conditions are equivalent\\
(i) $\O$ is weakly reversible,\\
(ii) $[\O ]_{\cong}$ is a convex set in the lattice $\la \Top _X ,\subset\ra$,\\
(iii) $[\O ]_\sim =[\O ]_{\cong}$.
\end{itemize}
\end{prop}
\dok
(a) If $\O _1, \O _2 \in [\O ]_\sim$ and $\O _1\subset \O '\subset \O _2 $, where $\O '\in \Top _X$,
then $\O _1 \sim \O _2$ and $\O _1\preccurlyeq  \O '\preccurlyeq \O _2 \preccurlyeq\O _1$,
which implies $\O '\sim \O _1$ and, hence, $\O '\in [\O _1 ]_\sim =[\O  ]_\sim$.

(b) Since $[\O ]_{\cong} \subset [\O ]_\sim$ by (a) we have $\Conv ([\O ]_{\cong})\subset [\O ]_\sim$.
Conversely, if $\O '\in [\O ]_\sim$,
then $\O \preccurlyeq \O' \preccurlyeq \O$.
Since $\O \preccurlyeq \O'$, by Fact \ref{T103} there is $\O _1 \in [\O ]_{\cong}$ such that $\O _1 \subset \O '$ and,
since $\O ' \preccurlyeq \O$, by Fact \ref{T103} again there is $\O _2 \in [\O ]_{\cong}$ such that $\O ' \subset \O _2$.
Thus $[\O ]_{\cong}\ni \O _1 \subset \O '\subset \O _2\in [\O ]_{\cong} $,
and, hence, $\O '\in \Conv ([\O ]_{\cong})$.

(c) The equivalence of (ii) and (iii) follows from (a) and (b), while (i) $\Leftrightarrow$ (iii) follows from Fact \ref{T103}.
\kdok
If $(X, \O _X)$ and $(Y, \O _Y)$ are homeomorphic spaces,
it is evident that the partial orders $\la [\O _X]_{\cong },\subset \ra$  and $\la [\O _Y]_{\cong },\subset \ra$ are isomorphic
and the same holds for the partial orders $\la [\O _X]_{\sim },\subset \ra$  and $\la [\O _Y]_{\sim },\subset \ra$.
Thus, in analogy with cardinal invariants,
the order types $\pi _{\cong}(\O):=\otp (\la [\O ]_{\cong },\subset \ra)$ and  $\pi _{\sim}(\O):=\otp (\la [\O ]_{\sim },\subset \ra)$
can be regarded as topological (order) invariants. So, a topology $\O$ is strongly reversible iff $\pi _{\cong}(\O)=\pi _{\sim}(\O)=1$;
more generally, $\O$ is reversible iff $\pi _{\cong}(\O)=\pi _{\sim}(\O)$ is an antichain. Otherwise we have
\begin{prop}\label{T106}
If $\O \in \Top _X$ is a non-reversible topology,
then each maximal chain $\L$ in $\la [\O ]_{\sim },\subset \ra$ is a Dedekind-complete linear order without end points.
\end{prop}
\dok
Suppose that $\O '$ is a largest element of $\L$.
By Proposition \ref{T104}(b) there are $\O _1,\O _2\in [\O ]_{\cong}$ such that $\O _1 \subset \O' \subset \O _2$
and, since $\O _2 \in [\O ]_{\sim}$, by the maximality of $\L$ we have $\O' = \O _2$.
Thus $\O '\in [\O ]_{\cong}$
and there is $f\in \Sym (X)$ such that $f [\O ]=\O '$.
Since the topology $\O$ is non-reversible
there is $\O _1 \in [\O ]_{\cong}$ such that $\O \varsubsetneq \O _1$
and, hence, $\O ' =f[\O ] \varsubsetneq f[\O _1]\in [\O ]_{\sim }$,
which is impossible by the maximality of $\L$.
Thus, $\L$ has no largest element and, similarly, $\L$ has no smallest element.

If $\la \L _0 ,\L _1\ra$ is a cut in $\la \L ,\subset\ra$, where $\L _0 ,\L _1 \neq \emptyset$,
then $\O ' :=\bigcap \L _1 \in \Top _X$
and by the convexity of $[\O ]_{\sim }$ we have $\O '\in [\O ]_{\sim }$.
Clearly $\la \L \cup \{ \O '\}, \subset \ra$ is a linear order
and by the maximality of $\L$ we have $\O '\in \L$;
thus, $\O '=\max \L _0$ or $\O '=\min \L _1$.
So $\L$ has no gaps; that is, it is Dedekind complete.
\kdok
\section{Weakly reversible non reversible spaces}
By Propositions \ref{T104} and \ref{T106}, if $\O \in \Top _X$ is a weakly reversible non-reversible topology,
then $[\O ]_{\sim }=[\O ]_{\cong }$ and each maximal chain in the poset $\la [\O ]_{\cong },\subset \ra$ is a Dedekind-complete linear order without end points.
Some topologies of that type and the corresponding maximal chains are described in the following example.
\begin{ex}\label{EX100}\rm
{\it $2^{\c}$ non-homeomorphic weakly reversible non-reversible topologies in $\Top _\c$.}
By Theorem 4.1 of \cite{JK}, for each infinite $1$-homogeneous\footnote{
A linear order  $\BL$ is called {\it $n$-homogeneous} (or $n$-transitive) iff for each $k\leq n$ and each $x_1,\dots ,x_k,$ $y_1,\dots, y_k\in L$,
such that $x_1<x_2<\dots <x_k$ and $y_1<y_2<\dots <y_k$,
there is an automorphism $f \in \Aut (\BL)$ such that $f(x_i)=y_i$, for all $i\leq k$.
$\BL$ is called {\it (ultra)homogeneous} if it is $n$-homogeneous, for each $n\in \N$, and this holds iff $\BL$ is $2$-homogeneous (see \cite{Ros}).}
Dedekind complete linear order $\BL $
there is a weakly reversible non-reversible topology $\O _\BL \in \Top _{|L|}$
such that
\begin{equation}\label{EQ104}\textstyle
\la [\O _\BL ]_{\cong}, \subset \ra =\la [\O _\BL ]_{\sim}, \subset \ra \cong \dot{\bigcup}_{2^{|L|}}\BL,
\end{equation}
(where $\dot{\bigcup}_{2^{|L|}}\BL$ denotes the disjoint union of $2^{|L|}$-many copies of $\BL$)
and, hence, each maximal chain in $[\O _\BL ]_{\cong}$ is isomorphic to $\BL$.
The topology $\O _\BL$ is constructed in the following way:
first we add a first element $z$ to $\BL$ and obtain the linear order $\BL _z=\{ z\}+\BL$,
then we take arbitrary $c\in L$ and define
\begin{equation}\label{EQ105}
\O _\BL :=\Big\{ [z,a):a\in L\Big\}\cup \Big\{ (z,b):b\in (z,c]\Big\}\cup \Big\{ \emptyset ,L_z\Big\}\in \Top _{L_z}.
\end{equation}
Concerning the diversity of $1$-homogeneous Dedekind complete linear orders
we note that the only scattered one is the integer line, $\Z$, and the non-scattered ones are dense and of size $\c$ \cite{Ohk}.
If $\BC$ is one of them and $(a,b)\cong\BC$, whenever $a<b\in C$,
then (see Fact 4.2 of \cite{JK}) the linear orders
$\BC$,
$R(\BC)  :=  \BC + (1+\BC)\o _1$,
$L(\BC)  :=  (\BC +1)\o _1^* +\BC$ and
$M(\BC)  :=  (\BC +1)\o _1^* +\BC+ (1+\BC)\o _1$
are $2$-homogeneous and Dedekind complete.
Of course, one of such $\BC$-s is the real line, $\R$,
and taking $\BL\in \{ \R ,R(\R), L(\R),M(\R)\}$ we have (\ref{EQ104}).

Further, taking for $\BC$ a homogeneous Aronszajn continuum (dense complete linear order) without its end points
(one was detected by Shelah \cite{She}, see \cite{Tod}, p.\ 260 for a construction),
for $\BL\in \{\BC ,R(\BC), L(\BC), M(\BC)\}$ we have (\ref{EQ104}).

Moreover,  Hart and van Mill \cite{Mil} constructed $2^{\c}$ non-isomorphic homogeneous continua, say $\BL _\a$, $\a <2^{\c}$.
Thus, if $\a <2^{\c}$ and $\BC _\a$ is the linear order $\BL _\a$ without its end points,
then for $\BL\in \{\BC_\a ,R(\BC_\a), L(\BC_\a), M(\BC_\a)\}$ we have (\ref{EQ104}).
If $\a <\b <2^{\c}$, then the posets $\dot{\bigcup}_{2^{\c}}\;\BC_\a$ and $\dot{\bigcup}_{2^\c}\;\BC_\b$ are not isomorphic
and, since $\pi _{\cong}(\O _{\BC_\a})\neq \pi _{\cong}(\O _{\BC_\b})$, the topologies $\O _{\BC_\a}$ and $\O _{\BC_\b}$ are not homeomorphic.
So we obtain $2^{\c}$ non-homeomorphic weakly reversible non-reversible topologies in $\Top _\c$.
\end{ex}
The weakly reversible non-reversible topologies defined by (\ref{EQ105}) are of simple form and irrelevant; 
they are not $T_1$, but they are first countable and, hence sequential\footnote{
A space is {\it sequential} iff a set $A$ is closed if $\lim \la x_n\ra \subset A$, for each sequence $\la x_n\ra$ in $A$.}.
Since the class of non-reversible spaces includes several prominent structures e.g.\ $\Q$, $\R \setminus \Q$ and normed spaces of infinite dimension (see \cite{Raj})
it is natural to ask are there some more relevant weakly reversible non-reversible spaces.
In the sequel we show that such examples can not be found in the class of sequential spaces in which each sequence has at most one limit\footnote{
These are exactly the spaces such that there is a partial function $\l $ from $X^\o$ to $X$ (a limit operator) 
such that $(X, \l)$ is an $\L ^*$-space and the sequential topology induced by $\l$ coincides with the original one (see \cite{Eng}, 1.7.18--1.7.20).};
thus, in particular, in the class of metrizable spaces.
We will use the following well-known consequence of Ramsey's theorem.
\begin{fac}\label{T101}
(a) If $f:\o \rightarrow \o$, then there is $H\in [\o]^\o$ such that the restriction $f\upharpoonright H$ is either a constant function
or a strictly increasing function.

(b) If $X\neq \emptyset$ and $f:\o \rightarrow X$,
then there is $H\in [\o]^\o$ such that the restriction $f\upharpoonright H$ is either a constant function
or an injection.
\end{fac}
\dok
(a) Defining $K_0:=\{ \{ k,l\}\in [\o ]^2: f(\min\{ k,l\})< f(\max\{ k,l\})\}$ and $K_1=[\o ]^2\setminus K_0$
by Ramsey's theorem there is $H\in [\o]^\o$, such that $[H]^2 \subset K_i$, for some $i<2$.
If $[H]^2 \subset K_0$, then for $k,l\in H$, where $k<l$ we have $f(k)<f(l)$ and $f\upharpoonright H$ is a strictly increasing function.
If $[H]^2 \subset K_1$ and $H=\{ k_i:i\in \o\}$ is an enumeration, where $k_0 <k_1,\dots$, then $f(k_0) \geq f(k_1),\dots$
and, hence, there are $i_0\in \o$ and $m\leq f(k_0)$ such that $f(k_i)=m$, for all $i\geq i_0$.

(b) Defining $K_0:=\{ \{ k,l\}\in [\o ]^2: f(\min\{ k,l\})\neq f(\max\{ k,l\})\}$ and $K_1=[\o ]^2\setminus K_0$
by Ramsey's theorem there is $H\in [\o]^\o$, such that $[H]^2 \subset K_i$, for some $i<2$.
If $[H]^2 \subset K_0$, then for $k,l\in H$, where $k\neq l$ we have $f(k)\neq f(l)$ and $f\upharpoonright H$ is an injection.
If $[H]^2 \subset K_1$, then for $k,l\in H$ we have $f(k)=f(l)$ and $f\upharpoonright H$ is a constant function.
\kdok
\begin{te}\label{T100}
If $\CX$ is a sequential space in which each sequence has at most one limit, then
\begin{equation}\label{EQ106}
\CX \mbox{ is weakly reversible }\Rightarrow \CX \mbox{ is reversible.}
\end{equation}
\end{te}
\dok
Let $\CX=(X, \O)$ and let $\F$ denote the corresponding family of closed sets. 

{\it Claim.} If $\la x_n :n\in \o\ra \rightarrow x$ and $x\not\in \{x_n :n\in \o \}$, then for each $M\in [\o ]^\o$ we have $F_M:=\{x_n :n\in M\} \cup \{ x\}\in \F$.

{\it Proof of Claim.}
Aiming for a contradiction, suppose that $F_M\not\in \F$.
Then there are a point $y\in X\setminus F_M$ and a sequence $\la b_k :k\in \o\ra$ in $F_M$ such that $\la b_k \ra\rightarrow y$.
Assuming that $|\{ k\in \o :b_k =x\}|=\o$,
the sequence  $\la b_k \ra$ would contain a constant subsequence $\la x,x,\dots\ra$ converging to $x$ and to $y$, which is impossible.
Thus, w.l.o.g.\ we suppose that $b_k =x_{n_k}$, for each $k\in \o$.
By Fact \ref{T101}(a) for the mapping $k\mapsto n_k$  there is a homogeneous set $H\in [\o]^\o$.
Assuming that there is $m\in \o$ such that $n_k=m$, that is, $b_k =x_m$, for all $k\in H$,
the sequence $\la b_k \ra$ would contain a constant subsequence $\la x_m,x_m,\dots\ra$ converging to $x_m$ and to $y$, which is impossible.
Thus the mapping $b\upharpoonright H$ is strictly increasing.
So if $H=\{ k_i:i\in \o\}$ is the enumeration, where $k_0 <k_1,\dots$,
then $n_{k_0}< n_{k_1}< \dots$,
which means that $\la b_{k_i} :i\in \o\ra =\la x_{n_{k_i}}:i\in \o\ra$
is a subsequence of the sequence $\la x_n\ra$ converging to $x$.
So, $\la b_{k_i} \ra$ converges to $x$ and to $y$, which is impossible.
Claim is proved.

We prove the contrapositive of (\ref{EQ106}).
Let $(X, \O)$ be a non-reversible space, $\O' \in [\O ]_{\cong}$, $\O \varsubsetneq \O '$ and $A\in \F ' \setminus \F$.
Since $(X, \O)$ is a sequential space and $A\not\in \F$
there are a point $x\in X\setminus A$ and a sequence $\la a_n :n\in \o\ra$ in $A$ such that $\la a_n\ra \rightarrow _\O x$.
By Fact \ref{T101}(b) we can suppose that the mapping $n\mapsto a_n$ is an injection.

By (a), for each $M\in [\o ]^\o$ we have $F_M:=\{ a_n:n\in M\}\cup \{ x\}\in \F \subset \F '$
and, since $A\in \F '$ we have $A_M:=F_M \cap A= \{ a_n:n\in M\}\in \F '\setminus \F$.
Let $\{ M_\a :\a <\c\} \subset [\o ]^\o$ be a mad family.
Then, since $X\setminus A_{M_\a}\in \O '$, for all $\a <\c$, we have
\begin{equation}\label{EQ100}
\O \subset \O ^* :=\O \Big(\O \cup \{ X\setminus A_{M_\a}:\a <\c\}\Big)\subset \O '.
\end{equation}
We will prove that $(X, \O ^*)$ is not a sequential space,
which will, by (\ref{EQ100}), imply that $(X, \O)$ is not weakly reversible.

First we show that in the space $(X, \O ^*)$ we have $x\in \overline{A_\o}$, where $A_\o =\{ a_n:n\in \o\}$.
By (\ref{EQ100}) a base for the topology $\O ^*$ consists of the sets of the form $O\setminus \bigcup _{\a\in K}A_{M_\a}$,
where $O\in \O$ and $K\in [\c ]^{<\o}$.
So, if $U\in \U (x)$ and $x\in O\setminus \bigcup _{\a\in K}A_{M_\a}\subset U$,
then for $\b \in \c \setminus K$ we have $S:=M_\b \setminus \bigcup _{\a\in K}M_\a\in [\o]^\o$
and, since the mapping $n\mapsto a_n$ is an injection, $A_S =\{ a_n :n\in S\}\subset X\setminus \bigcup _{\a\in K}A_{M_\a}$.
Now $\la a_n :n\in S\ra$ is a subsequence of the sequence $\la a_n :n\in \o\ra$
and, hence, $\la a_n :n\in S\ra\rightarrow _\O x$;
so, since $x\in O \in \O$,
there is $n_0\in \o$ such that for each $n\in S\setminus n_0$ we have $a_n\in O\setminus \bigcup _{\a\in K}A_{M_\a}\subset U$.
Thus $x\in \overline{A_\o}\setminus A_\o$ indeed.

Second, we show that the set $A_\o$ is sequentially closed in the space $(X, \O ^*)$.
On the contrary, suppose that there are a point $y\in X\setminus A_\o$ and a sequence $\la b_k :k\in \o\ra$ in $A_\o$ such that $\la b_k \ra \rightarrow _{\O^*}y$.
Then $\la b_k \ra \rightarrow _\O y$ and, hence, $y\in \overline{A_\o}^\O$.
By (a) we have $F_\o =A_\o \cup \{ x\}\in \F$, which gives $\overline{A_\o}^\O=A_\o \cup \{ x\}$
and, consequently, $y=x$ and $\la b_k \ra \rightarrow _{\O^*}x$.
So, $\la b_k :k\in \o\ra=\la a_{n_k} :k\in \o\ra$
and, since $\la b_k :k\in \o\ra$ has no constant subsequence,
by Fact \ref{T101}(a) there is a set $H=\{ k_i:i\in \o\}\in [\o ]^\o$, where $k_0 <k_1 <\dots$,
such that $n_{k_0} <n_{k_1} <\dots$.
Thus we have $\la a_{n_{k_i}} :i\in \o\ra\rightarrow _{\O^*}x$.
Since $M:=\{ n_{k_i} :i\in \o\}\in [\o ]^\o$
there is $\a <\c$ such that $|M\cap M_\a|=\o$
and, since the mapping $n\mapsto a_n$ is an injection,
\begin{equation}\label{EQ101}
|A_M\cap A_{M_\a}|=|\{ a_{n_{k_i}}: i\in \o \land n_{k_i}\in M_\a \}|=\o .
\end{equation}
Now we have $U:=X\setminus  A_{M_\a}\in \U _{\O ^*}(x)$
and, hence, there is $i_0\in \o$ such that for each $i\geq i_0$ we have $a_{n_{k_i}}\not\in A_{M_\a}$
which contradicts (\ref{EQ101}).
\kdok

\section{Strongly reversible spaces}
We recall that a space $\CX =(X,\O)$ is called strongly reversible
iff each bijection $f:\CX \rightarrow \CX$ is a homeomorphism;
that is, $f[\O]=\O$, for each $f\in \Sym (X)$. Thus $[\O ]_{\cong}=\{ f[\O] :f\in \Sym (X)\}=\{ \O\}=[\O ]_{\sim}$.
Here we show that the class of strongly reversible topologies on a set $X$ is small:
it contains discrete and antidiscrete topologies and generalizations of the cofinite topology.

If $X$ is a set of size $\k$ and $\l,\mu \leq \k$ are cardinals,
we use standard notation $[X]^\mu:=\{ A\subset X: |A|=\mu\}$, $[X]^{<\mu}:=\{ A\subset X: |A|<\mu\}$ etc.
If, in addition, $\l +\mu =\k$, we define $[X]^{\l |\mu}:=\{ A\subset X: |A|=\l \land |X\setminus A|=\mu\}$,
and it is evident that for $A,B\in [X]^{\l |\mu}$ there is $f\in \Sym (X)$ such that $f[A]=B$.
For $\C \subset P(X)$ and $f\in \Sym (X)$ let $f[\C]=\{ f[C]:C\in \C\}$.
\begin{te}\label{T102}
A topology $\O\in \Top _X$ is strongly reversible iff it is discrete or antidiscrete or $|X|\geq \o$ and
$\O =\{ X\setminus F : F\in [X]^{<\l}\}\cup \{ \emptyset\}$, for some cardinal $\l \in [\o ,|X|]$.
\end{te}
\dok
Let $|X|=\k$ and let $f[\O ]=\O$, for all $f\in \Sym (X)$.

If there exists $F\in (\F \setminus \{ X \})\cap [X]^\k$,
then $0<\mu := |X\setminus F|\leq \k \geq \o$
and $X\setminus F \in \O \cap [X]^{\mu |\k}$.
So, for each $O'\in [X]^{\mu |\k}$ there is $f\in \Sym (X)$ such that $O'=f[X\setminus F]$;
thus, $[X]^{\mu |\k}\subset \O$.
Let $x\in X$ and let $\nu :=\mu -1$, if $\mu <\o$; and $\nu :=\mu $, if $\mu \geq\o$.
Since $\nu \leq \k \geq \o$ there are disjoint sets $A_0,A_1\in [X\setminus \{ x\}]^\nu$
and defining $O_i=\{ x\}\cup A_i$, for $i<2$,
since $|\nu +1|=\mu$, we have $O_i \in [X]^{\mu |\k}\subset \O$, for $i<2$,
and, hence, $\{ x\}=O_0 \cap O_1 \in \O$.
Thus $\{ x\} \in \O$, for all $x\in X$,
and, hence, $\O =P(X)$.

Otherwise we have $\F \setminus \{ X \} \subset [X]^{<\k}$ and, hence,
\begin{equation}\label{EQ102}
 \l :=\min\{ \mu : \F \setminus \{ X \} \subset [X]^{<\mu }\}\leq \k \;\;\mbox{ and }\;\;\F \setminus \{ X \} \subset [X]^{<\l}.
\end{equation}
So, if $\l=1$, then $\F =\O =\{ \emptyset ,X\}$ and we are done.

If $\l >1$, we prove that  $[X]^{<\l}\subset \F \setminus \{ X \}$,
which will, together with (\ref{EQ102}), give $\O =\{ X\setminus F : F\in [X]^{<\l}\}\cup \{ \emptyset\}$;
thus, we show that
\begin{equation}\label{EQ103}
\forall \mu <\l \;\;  [X]^{\mu}\subset \F.
\end{equation}
So, if $1\leq \mu <\l$,
then by (\ref{EQ102}) there is $F\in \F \setminus \{ X \}$ such that  $\mu \leq |F|=:\theta <\l \leq \k$.
Since $\theta <\k$ the equation $\theta +\nu =\k$ has a unique solution $\nu$
(if $\k <\o$ that is evident; otherwise we have $\nu =\k$)
and, consequently, $[X]^\theta =[X]^{\theta |\nu} \ni F$.
Now, for each $G\in [X]^{\theta }$ there is $f\in \Sym (X)$ such that $G=f[F]\in \F$
and, thus, $[X]^{\theta }\subset \F$.
So, if $\theta =\mu$, we are done.
If $\theta >\mu$ and $B\in [X]^\mu$,
then there are $F_0,F_1 \in [X]^{\theta }$ such that $B=F_0 \cap F_1$
and, hence, $B\in \F$.
Thus (\ref{EQ103}) is true.

Finally, since $\l >1$ by (\ref{EQ103}) we have $[X]^1 \subset \F$
and, since the family $\F$ is closed under finite unions, $[X]^{<\o}\subset \F$.
So if $\k <\o$, then $\O =P(X)$
and if $\k \geq \o$, then $\l \geq \o$ and $\O =\{ X\setminus F : F\in [X]^{<\l}\}\cup \{ \emptyset\}$.
\kdok

\paragraph{Acknowledgement.}
This research was supported by the Science Fund of the Republic of Serbia,
Program IDEAS, Grant No.\ 7750027:
{\it Set-theoretic, model-theoretic and Ramsey-theoretic
phenomena in mathematical structures: similarity and diversity}--SMART.

\end{document}